\documentclass[12pt]{amsart}
\usepackage{amsmath}
\usepackage{amsfonts,amssymb}
\usepackage{amsthm}
\usepackage[matrix,arrow,curve]{xy}
\usepackage{tikz}
\usetikzlibrary{positioning}
\tikzset{>=stealth}
\usetikzlibrary{arrows}
\usetikzlibrary{calc}
\usetikzlibrary{decorations.markings}
\usetikzlibrary{matrix}
\numberwithin{equation}{section}
\theoremstyle{plain}
\newtheorem{theorem}{Theorem}[section]

\newtheorem{conjecture}[theorem]{Conjecture}
\theoremstyle{definition}

\newtheorem{remark}[theorem]{Remark}

\usepackage{array}

\begin{document}
\title[about the structure]{On the moduli space of Higgs pairs with nilpotent residues.}

\author{Denis Degtyarev}

\address{
National Research University Higher School of Economics, Moscow, RUSSIA}
\email{degtyarev.d.o@gmail.com}

\begin{abstract}We specify the conjecture about the structure of
the Coulomb branch of N=2 supersymmetric quantum field theories in four dimensions
\end{abstract}

\maketitle

\tableofcontents

\section{Introduction}

Let us given a N=2 supersymmetric quantum field theory (4d N=2 SUSSY QFT) $Q$. Physical arguments imply that we can associate to $Q$ an algebraic integrable system $DW(Q)$ which is in fact a generalized Hitchin system \cite{Don1}. More precisely, let $\Sigma$ be a smooth algebraic curve. In this paper we study the moduli space of meromorphic $G$-Higgs pairs on $\Sigma$. The moduli space of such bundles has a canonical Poisson structure. Symplectic leaves of the Poisson structure are the main objects of our interest. 

We start by recalling basic constructions of the moduli space of Higgs bundles and Poisson structure in section 2, following \cite{Mar94}. A point in a symplectic leaf of the Poisson structure is an isomorphism class of pair $(E,\varphi)$, where $E$ is the principal $G$-bundle on $\Sigma$ and $\varphi \in H^{0}(\Sigma, \mathrm{ad}(E)\otimes K(D))$ is the meromorphic section with poles in the finite set of points $D$ on $\Sigma$. Moreover, residues of meromorphic sections in a symplectic leaf belong to fixed coadjoint orbits of Lie algebra $\mathfrak{g}$.   In section 3 we concentrate on those symplectic leaves for which meromorphic sections of $G$-Higgs bundles have nilpotent residues.

One of the conjectures in \cite{Tachikawa1} essentially states that the Hitchin map from $DW(Q)$ factors through an affine space $\mathcal{M}_{Coulomb}(Q)$. In section 4, we propose a further detailed version of the conjecture about the structure of the Coulomb branch of the theory $Q$ suggested in \cite{Tachikawa1}. According to our conjecture, we identify the Coulomb branch of the theory $Q$ with the affinization of the suitably defined finite cover of the symplectic leaf of the moduli space of meromorphic $G$-Higgs bundles.

\section{Preliminaries}

\subsection{Generalities on moduli space of Higgs bundles.}

Let $\Sigma$ be a smooth complete algebraic curve of genus $g$, $K$ its canonical line bundle. Let $G$  be a complex semisimple Lie group with the Lie algebra $\mathfrak{g}$.

\subsubsection{Meromorphic Higgs pairs.}
\textit{A holomorphic $K$-valued $G$-Higgs pair} on $\Sigma$ is a pair $(E,\varphi)$ consisting of principal holomorphic $G$ bundle $E$ and a holomorphic section of the twist $\mathrm{ad}E\otimes K$ of its adjoint Lie algebra bundle by canonical line bundle, i.e. $\varphi \in H^0(\Sigma, \mathrm{ad}E\otimes K)$. Let $D$ be a closed zero dimensional subscheme of $\Sigma$. We will denote also by $D$ the corresponding effective divisor $D:=p_1+\dotsc+p_m$  on $\Sigma$, which consist of distinct points. The pair $(E,\varphi)$ is a meromorphic Higgs pair if $\varphi \in H^0(\Sigma,\mathrm{ad}(E)\otimes K(D))$. For simplicity, we shall denote $K(D)$ by $L$ .

\subsubsection{Stable Higgs bundles.}  
\textit{The Higgs pair is (semi-) stable}, if for any $\varphi$-invariant reduction $\sigma: X \to E/P$ of $E$ to a maximal parabolic subgroup $P \subset G $, $\mathrm{deg}\sigma^{*}T_{G/P} > 0$ (resp. $\geqslant 0$), for more details cf. \cite{Ram75}. Here $T_{G/P}$ is a vertical bundle of
$E/P \to X$. 

According to \cite{Ram75}, there is an analytic coarse moduli space $\mathcal{M}(G,D) = \coprod\limits_{c\in \pi_1(G)}\mathcal{M}(G,D,c)$ of
semi-stable pairs $(E,\varphi)$. Connected components of $\mathcal{M}(D)$ are labelled by the topological type of $G$ bundle $E$. We will denote the component of topologically trivial $G$-bundles also by $\mathcal{M}(D)$ for short.

\subsection{Poisson structure via levels.} 
In the classical setting  of \cite{Hit1} the moduli space of holomorphic stable Higgs pairs is a partial compactification of the cotangent bundle of the moduli space of stable vector bundles. Therefore, it has a canonical symplectic structure. The moduli space $\mathcal{M}(D)$ of stable meromorphic Higgs pairs has a natural Poisson structure, see \cite{Bot1},\cite{Mar94}. We describe a Poisson structure on $\mathcal{M}(D)$ in two steps in much the same way as in the \cite{Mar94} but in a slightly different generality. 

\begin{enumerate}
\item As a preliminary step we realize a dense open subset of $\mathcal{M}(D)$ as an orbit space of a Poisson action of a group $G_D$ on the cotangent bundle of the moduli space $G$-bundles with level structures.

\item Next we exhibit a two-vector on the smooth locus of $\mathcal{M}(D)$ through the use of the cohomological construction.
\end{enumerate}
\subsubsection{D-level structure.} 
\textit{A $D$-level structure} on a $G$-bundle $E$ is an isomorphism $\eta: E|_D \to G\times D$. Let $\mathcal{U}(\Sigma,G,c,D) = \coprod\limits_{c\in \pi_1(G)}\mathcal{U}(G,D,c)$ be the moduli space of $ G$-principal bundles with a level structure. We will denote connected component with $c=0$ by $\mathcal{U}(D)$ for short. A point in $\mathcal{U}(D)$ is an isomorphism class of pair $(E,\eta)$. The group $\widetilde{G}_D$ of maps from $D$ to $G$ acts naturally on $\mathcal{U}(D)$ by changing the framing data $\eta$: an element $\widetilde{f} \in \widetilde{G}_D$ sends $[(E,\eta)] \to [(E, \widetilde{f} \circ \eta)]$ where $\widetilde{f}$ lifts $f$, $[\cdot]$ denotes the isomorphism class of a pair. The action of $G_D$ factors through its quotient $G_D:=\widetilde{G}_D/Z(G)$ by the diagonal embedding of the center of $G$. We will refer to $G_D$ as a level group.

\subsubsection{The cotangent bundle of $\mathcal{U}(D)$.} Deformation theory arguments give identification $T^{\vee}_{[(E,\eta)]}\mathcal{U}(D)\cong H^1(\Sigma, \mathrm{ad}(E)(-D))$. Using Serre's duality we obtain natural isomorphism:
\begin{equation}
T^{\vee}_{[(E,\varphi)]}\mathcal{U(D)}\simeq H^0(\Sigma, \mathrm{ad}E\otimes K(D)),
\end{equation}
a point in $T^{\vee}\mathcal{U}(D)$ is an isomorphism class of triple $[(E,\varphi,\eta)]$ where $(P,\varphi)$ is a Higgs pair.

The action of the level group $G_D$ on $\mathcal{U}(D)$ lifts naturally to an action on $T^{\vee}\mathcal{U}(D)$:
\begin{equation}
[(E,\varphi,\eta)] \mapsto [(E,\varphi, \widetilde{f} \circ \eta)],
\end{equation}

where $\widetilde{f} \in \widetilde{G}_D$ lifts $f$. The lifted action exhibits several remarkable properties:
\begin{enumerate}
\item $T^{\vee}\mathcal{U}(D)$ carries a Poisson action $G_D$ with respect to the canonical symplectic form on the cotangent bundle $T^{\vee}\mathcal{U}(D)$.  
 
\item The moment map $\mu:T^{\vee}\mathcal{U}(D) \to \mathfrak{g}_D^{\vee}$ is given by
\begin{equation}\label{moment map}
\mu (E,\varphi,\eta): A \mapsto \mathrm{Res}\,\mathrm{Tr}(A\cdot\varphi^{\eta}), 
\end{equation}
where $\varphi^{\eta}:=\eta\circ \varphi \circ \eta^{-1} \in H^{0}(\Sigma, \mathfrak{g}\otimes K(D)|_D)$, $A \in \mathfrak{g}_D$ and the residue map $\mathrm{Res}:H^{0}(K(D)\otimes \mathcal{O}_D) \to H^{1}(\Sigma, K )\simeq\mathbb{C}$.
\item $G_D$ acts freely on the open subset $(T^{\vee}\mathcal{U}(D))^{\circ}$ parametrizing isomorphism classes of triples $[(E,\varphi,\eta)]$ where $(E,\eta)$ is a stable pair of $G$-bundle with $D$-level structure and $(E,\varphi)$ is a stable Higgs pair.
\end{enumerate}

The third property of the lifted action makes $T^{\vee}\mathcal{U}$
into principal $G_D$-bundle over an open subset $\mathcal{M}^{\circ}(D)$ of stable Higgs pairs $\mathcal{M}^s(D)$. Hence the canonical symplectic structure induces a Poisson structure on $\mathcal{M}^{\circ} \simeq T^{\vee}\mathcal{U}^{\circ}(D)/G_D$. 

A canonical global formula for this Poisson structure is obtained via cohomological  techniques. The tangent space to $\mathcal{M}(D)$ at a smooth point is identified with the first hyper-cohomology of the complex $\mathcal{K}:=\mathrm{ad}(E)\xrightarrow{[\varphi,\cdot]}\mathrm{ad}(E)\otimes L$:
\begin{equation}
T_{[(E,\varphi)]}\mathcal{M}(D) \cong \mathbb{H}^1( \mathcal{K})
\end{equation}
The dual of the complex $\mathcal{K}$ is the complex

\begin{equation}
\mathcal{K}^{\vee}:=(\mathrm{ad}(E)\otimes L)^{\vee}\otimes K \xrightarrow{[\varphi,\cdot]^{\vee}}(\mathrm{ad}(E))^{\vee}\otimes K
\end{equation}

Grothendieck duality gives a natural isomorphism $\mathbb{H}^1(\mathcal{K}) \cong \mathbb{H}^1( \mathcal{K}^{\vee})^{\vee}$. For $\mathcal{K}$ we obtain the canonical isomorphism of complexes $\mathcal{K}^{\vee}\cong\mathcal{{K}}\otimes \mathcal{O}(-D)$. Composing with embedding $i: \mathcal{O}(-D)\to \mathcal{O}$ we get the morphism of complexes:
\begin{equation}
I: \mathcal{K}^{\vee} \to \mathcal{K},
\end{equation}
and using the duality we obtain the map
\begin{equation}
I: \mathbb{H}^1(\mathcal{K})^{\vee}\cong \mathbb{H}^1(\mathcal{K}^{\vee})\to \mathbb{H}^1(\mathcal{K}).
\end{equation}
The map $I$ gives the element $\Omega_D\in\otimes^2\mathbb{H}^1(\mathcal{K})$ which is skew-symmetric. In such a way we have the global two-vector  $\Omega_D$ on the nonsingular part $\mathcal{M}^{ns}(D)$ of the moduli space $\mathcal{M}(D)$. It turns out that restriction of $\Omega_D$ to $\mathcal{M}^{\circ}(D)$ is induced by the above Poisson structure on the quotient $T^{\vee}\mathcal{U}^{\circ}(D)/G_D$.
The two-vector is then automatically Poisson. 
\section{Integrable systems on symplectic leaves and their covers.}
\subsection{Marsden-Weinstein reduction}
The moment map $\ref{moment map}$ determines the foliation by symplectic leaves on $\mathcal{M}^s(D)$. Let $\mu([(E,\varphi,\eta)]) = e$, where $[(E,\varphi,\eta)]\in T^{\vee}\mathcal{U}(D)$, let $\mathcal{O}_e \subset \mathfrak{g}^{\vee}_D$ be the coadjoint orbit through $e$ and let $G_D^e$ be the stabilizer of $e$. Then the symplectic leaf $\mathcal{M}(\mathcal{O}_e)$ through $[(E,\varphi,\eta)]$  is the connected component of  
\begin{equation}\label{Marsden iso}
\mu^{-1}(\mathcal{O}_e)/G_D \cong \mu^{-1}(e)/G_D^e,
\end{equation}
see \cite{MarWei} for details . 
Markman made remarkable observation in \cite{Mar94} that the moduli space of stable meromorphic Higgs pairs $\mathcal{M}^s(D)$ is a completely integrable system. The lagrangian fibration of $\mathcal{M}^s(D)$ is given by the Hitchin map $h: \mathcal{M}^s(D) \to \bigoplus\limits_{i=1} H^{0}(\Sigma, L^{d_i})$ where $d_i$ is the degree of the generator $h_i$ of the algebra $\mathbb{C}[h_1,\dotsc, h_r]$ of invariant polynomials on $\mathfrak{g}$ with $r=\mathrm{rk}\,\mathfrak{g}$. The symplectic leaf $\mathcal{M}(\mathcal{O}_e)$ under the simplifying technical assumptions on the orbit $\mathcal{O}_e$ admits analogous construction, see \cite{Mar99}. In particular, the restriction of the Hitchin map $h$ to the symplectic leaf $\mathcal{M}(\mathcal{O}_e)$ is a Lagrangian fibration. 
\begin{remark}
An element $[(E,\varphi)] \in \mathcal{M}^s(D)$ belongs to the symplectic leaf $\mathcal{M}(\mathcal{O}_e)$ if  $\mathrm{Res}_{p_j}\varphi \in \mathcal{O}_{e_j}$ at each $p_j$, $j=1,\dotsc, m$.
\end{remark}

The cohomological identification of the Possion structure affords an opportunity to determine its rank on the symplectic leaf $\mathcal{M}(\mathcal{O}_e)$ (\cite{Mar94}). The rank of $\Omega_D$ at $[(E,\varphi)] \in \mathcal{M}(D)$ is equal to 
\begin{equation}
\mathrm{rk}\,\Omega_D=2(g-1)\mathrm{dim}\,G + \sum\limits_{j=1}^m \mathrm{dim}\,\mathcal{O}_{e_j},
\end{equation}
as a consequence dimension of the generic fiber equal to $\frac12 \Omega_D$.

\subsection{Covers of symplectic leaves.} Let $\mathfrak{g}$ be a semisimple algebra of classical type. $G$ will denote the Lie group of adjoint type.
\subsubsection{Nilpotent orbits.}From now on we assume that  coadjoint orbit $\mathcal{O}_e$ is nilpotent, i.e. each component $\mathcal{O}_{e_j} \in \mathfrak{g}$ in $\mathcal{O}_e \cong \mathcal{O}_{e_1}\times\dotsc\times \mathcal{O}_{e_m}$ is an orbit of the nilpotent element $e_j \in \mathfrak{g}$. In this situation the component group $A(\mathcal{O}_{e_j}):= G^{e_j}/(G^{e_j})^0$ is either trivial or a finite product of $\mathbb{Z}/2$, see \cite{ColMC}. The precise structure of $A(\mathcal{O}_{e_j})$ is completely determined by the partition $\mathrm{d}^j:=[\mathrm{d}_1^j,\dotsc, \mathrm{d}_{N_j}^j]$ of the nilpotent orbit $\mathcal{O}_{e_j} = \mathcal{O}_{\mathrm{d}^j}$. 

\subsubsection{Covers.}
In view of the right hand side of the isomorphism $\ref{Marsden iso}$, it is natural to consider a covering space for $\mathcal{M}(D)$. We will denote by $\widetilde{\mathcal{M}}(\mathcal{O}_e)$ the quotient  $\mu^{-1}(e)/(G^{e}_D)^0$. Let $p: \widetilde{\mathcal{M}}(\mathcal{O}_e) \to \mathcal{M}(\mathcal{O}_e)$
be the projection map.

\section{Coulomb branch}
Now deal with the restriction of the Hitchin map to the symplectic leaf $\mathcal{M}(\mathcal{O}_e)$ for the nilpotent orbit $\mathcal{O}_e$:
\begin{equation}
 h = (h_1,\dotsc,h_r):\mathcal{M}^s_{\mathcal{O}}\to B:=\bigoplus\limits_{i=1}^{r}H^{0}\left(\Sigma, K^{d_i}+ \sum\limits_{j=1}^{m}\nu_{d_i}(p_j)p_j\right)
\end{equation}
where $\nu_{d_i}(p_j)$ is the pole order of the $d_i$-differential at the point $p_j$. Let $X^{\mathrm{aff}}:=\mathrm{Spec}(X,\mathcal{O}_X)$ denote the affinization of an algebraic variety, let $\mathbb{C}[X]$ denote the ring of regular functions on $X$. We propose the following.

\begin{conjecture}
 The canonical morphism $f:\mathcal{M}(\mathcal{O}_e)\to \mathcal{M}^{\mathrm{aff}}(\mathcal{O}_e)$ can be included into the commutative diagram:
 \begin{center}
 \begin{tikzpicture}
   \matrix (m) [matrix of math nodes,row sep=3em,column sep=4em,minimum width=2em]
   {
      {\widetilde{\mathcal{M}}}_{\mathcal{O}_e} &  \mathcal{M}_{\mathcal{O}_e}^{\mathrm{aff}} \\
      \mathcal{M}_{\mathcal{O}_e} &  B\\};
   \path[-stealth]
     (m-1-1) edge node [above] {$\widetilde{h}$} (m-1-2)
     (m-2-1) edge node [above] {$h$} (m-2-2)
     (m-1-1) edge node [left]{$p$} (m-2-1)
     (m-1-2) edge node [right] {$f$} (m-2-2);

 \end{tikzpicture}
 \end{center}
 where the finite morphism $f$ is induced by the ring morphism  $u: \mathbb{C}[\mathcal{M}_{\mathcal{O}_e}] \to \mathbb{C}[\widetilde{\mathcal{M}}_{\mathcal{O}_e}]$.
\end{conjecture}
\begin{remark}The main value of the conjecture  is that we obtain the isomorphism $\mathcal{M}_{Coulomb} \cong \mathcal{M}^{\mathrm{aff}}_{\mathcal{O}_e}$.  It is remarkable that $\mathcal{M}^{\mathrm{aff}}_{\mathcal{O}_e}$ has expected dimension. In accordance with the above conjecture and the expression for the rank of $\Omega_D$  we have $\mathrm{dim}\, \mathcal{M}^{\mathrm{aff}}_{\mathcal{O}_e} = \frac12\mathrm{rk} \Omega_D  = (g-1) \mathrm{dim}G + \sum\limits_{j=1}^m \frac12\mathrm{dim}\,\mathcal{O}_{e_j}$. 
\end{remark}

\end{document}